\documentclass[11pt]{amsart}

\usepackage{graphics}
\usepackage{epsfig}
\usepackage{amsmath}
\usepackage{amscd}
\usepackage{amsfonts,latexsym,amssymb}

\newcommand{\SN}{{\mathcal{N}}}
\newcommand{\SO}{{\mathcal{O}}}

\newcommand{\PP}{\mathbb{P}}
\newcommand{\ZZ}{\mathbb{Z}}

\newcommand{\RR}{\mathbb{R}}

\newcommand{\isom}{\cong}
\newcommand{\Ext}{\operatorname{Ext}}

\newcommand{\Hom}{\operatorname{Hom}}

\newcommand{\NS}{\operatorname{NS}}

\newcommand{\id}{\operatorname{I}}

\newcommand{\inj}{\hookrightarrow}

\newcommand{\rk}{\operatorname{rk}}

\newcommand{\wt}{\widetilde}
\newcommand{\del}{{\partial}}
\newcommand{\delb}{{\overline\partial}}

\newcommand{\Res}{\operatorname{Res}}
\newcommand{\ch}{\operatorname{ch}}
\newcommand{\tr}{\operatorname{tr}}
\newcommand{\pardeg}{\operatorname{par-deg}}

\newcommand{\diag}{\operatorname{diag}}
\newcommand{\minus}{{-}}
\newcommand{\codim}{\operatorname{codim}}

\newtheorem{proposition}{Proposition}[section]
\newtheorem{theorem}[proposition]{Theorem}

\newtheorem{lemma}[proposition]{Lemma}

\newtheorem{corollary}[proposition]{Corollary}

\begin{document}

\title[Parabolic bundles]{Parabolic bundles 
and representations of the fundamental group}

\author[T. G\'omez, T. R. Ramadas]{Tom\'as L. G\'omez and T. R. Ramadas}

\address{School of Mathematics, 
Tata Institute of Fundamental
Research, Homi Bhabha road, Mumbai 400 005, India}
\email{tomas@math.tifr.res.in, ramadas@math.tifr.res.in}

\date{1 December 1999}

\begin{abstract}
Let $X$ be a smooth complex projective variety with Neron-Severi group
isomorphic to 
$\ZZ$, and $D$ an irreducible divisor with
normal crossing singularities. Assume $1<r\leq 3$. 
We prove that if 
$\pi_1(X)$ doesn't have irreducible $PU(r)$ representations, then 
$\pi_1(X\minus D)$ doesn't have irreducible $U(r)$ representations. 
The proof
uses the non-existence of certain stable parabolic bundles.
We also obtain a similar result for $GL(2)$ when $D$
is smooth and $X$ is a complex surface.

\end{abstract}

\maketitle

\section*{Introduction}

Let $X$ be a smooth complex projective variety of dimension $n$. 
Let $D
\subset X$ be an irreducible divisor with normal crossing
singularities. One would like to relate the fundamental groups of
$X$ and $X\minus D$. There is a short exact sequence
$$
1 \to N \to \pi_1(X\minus D) \to \pi_1(X) \to 1.
$$ 
Fix once and for all an element $\lambda$ of $\pi_1(X\minus D)$
going once around $D$.
The kernel $N$ is generated by the set 
$$
\{ a \lambda^{\sigma} a^{-1}: \sigma=\pm 1,\; a\in\pi_1(X\minus
D)
\}.
$$
The most definitive results on 
the fundamental group of $X \minus D$
are due to Nori. If $X$ is a surface 
\cite[prop. 3.27]{No} implies that if $D^2>2 r(D)$ (where $r(D)$
is number of nodes of $D$) then $N$ is a finitely generated abelian
group, and its centralizer is a subgroup of finite index. In
particular, his result implies Zariski's conjecture: if $X=\PP^2$,
then $\pi_1(\PP^2\minus D)$ is abelian (since in this case 
$D^2>2 r(D)$ is automatically satisfied).

We will make the
following assumption on the Neron-Severi group $\NS(X)$ and the rank
$r$ of the representations: 
\begin{equation}
\label{assumption}
\NS(X)\isom \ZZ L \quad \text{and} \quad 1<r\leq 3.
\end{equation}
The main result of this paper is that if $\pi_1(X)$ has no irreducible $PU(r)$ 
representations, then 
$\pi_1(X\minus D)$ has no irreducible $U(r)$ representations.

The motivation for this kind of result is the following.
If $\rho:\pi_1(X\minus D) \to U(r)$ is a representation such that 
$\rho(\lambda)$ is a multiple of the identity, then
$\rho$ descends to give a representation 
$\overline\rho:\pi_1(X)\to PU(r)$. If $\overline\rho$ is reducible,
then also $\rho$ is reducible. Of course this argument doesn't work if
$\rho(\lambda)$ is not a multiple of the identity, 
but using parabolic bundles (and assuming (1)) we show that this
cannot happen (corollary \ref{corollary}).

We also give a similar result 
for $GL(2)$
representations (theorem \ref{irred2}). Here we have to assume that
$X$ is a surface and $D$ is smooth, since the correspondence 
between
representations and parabolic Higgs bundles is only known under those
conditions.

\section{Unitary representations of fundamental groups}
\label{unitary}

Let $(X,H)$ be a polarized smooth projective variety. 
All
degrees and stability will be with respect to the polarization $H$,
unless otherwise stated.
Let $D
\subset X$ be an irreducible divisor with normal crossing
singularities.
Let 
$$
\pi:\wt D \to D \subset X
$$
be the composition of the normalization of $D$ with the inclusion in $X$.
Let $\rho:\pi_1(X \minus D)\to U(r)$ be a representation of the
fundamental group of the complement. This gives a local system on
$X \minus D$. Let $(E,\nabla)$ be the Deligne extension
(\cite{D}, \cite{Ka}),
i.e. $E$ is a holomorphic vector bundle on $X$ (with
$\delb_E=\nabla^{0,1}$
) and $\nabla:\!E \to E \otimes \Omega_X
\langle \log D\rangle  $ is a holomorphic logarithmic connection. 
The restriction
of the residue $\Res (\nabla): E|_{\wt D} \to  E|_{\wt D}$ to a point
in the normalization $\wt D$ of $D$ is 
(up to conjugation) $\Gamma$,
where $\exp (-2\pi i \Gamma)$ is the holonomy of $\rho$ on a small
loop around $D$, and the eigenvalues of $\Gamma$ satisfy $0 \leq
\alpha_1 \leq \alpha_2 \leq \ldots\leq \alpha_r < 1$. 

Recall that a parabolic structure on a vector bundle $E$ or rank $r$ 
is a filtration of
$E|_{\wt D}$ by holomorphic subbundles
$$
E|_{\wt D}=\pi^* E=F_1 \varsupsetneq F_2 \varsupsetneq \ldots 
\varsupsetneq F_p \varsupsetneq F_{p+1}=0
$$
with weights
$$
0 \leq \mu_1 < \mu_2 <\ldots <\mu_p <1
$$
where $\mu_j$ are the different eigenvalues (without repetitions) of
$\Gamma$. We say that the parabolic structure is trivial if $p=1$
(i.e. all weights $\alpha_i$ are equal).

If $E$ is a torsion-free sheaf, then $E|_{\wt D}$ won't be locally
free in general. In this case we only define the filtration in a
Zariski open subset of $\wt D$, in which $E|_{\wt D}$ is locally free. This
will be sufficient for our purposes (for a more general definition and
for the notion of morphism between parabolic sheaves, see
\cite{M-Y}).
 
The parabolic degree of $E$ is defined as
$$
\pardeg E=\deg E +\Big( \sum_{i=1}^r \alpha_i \Big) \deg D.
$$ 

We say that a parabolic bundle $E$ is parabolic stable
(resp. semistable) if for any saturated torsion-free parabolic 
subsheaf $E'$ or rank $r'$,
$$
\frac{\pardeg E'}{r'} < \frac{\pardeg E}{r} ,\quad (\text{resp.} \;
\leq).
$$
And in general, all notions related to Mumford stability (slope,
polystability, etc...) have a corresponding parabolic notion, changing
the usual degree with the parabolic degree.

\begin{proposition} 
\label{biquard}
Let  $\rho:\pi_1(X\minus D) \to
U(r)$ be a representation. Then there is an associated polystable parabolic
bundle $E$ of rank $r$. If the representation 
$\rho$ is irreducible then $E$ is parabolic stable. 

The Chern characters of $E$ are given by \cite[cor (B.3)]{E-V}
\begin{equation}
\label{chernchar}
\ch_1(E)=-\tr(\Gamma)[D],\qquad 
\ch_2(E)=\frac{1}{2} \tr(\Gamma^2)[D]^2.
\end{equation}
Equivalently, we have the Chern classes
\begin{equation}
\label{chern}
c_1(E)=-\tr(\Gamma)[D],\qquad 
c_2(E)=\frac{1}{2} \Big( (\tr\ \Gamma)^2-\tr\ (\Gamma^2) \Big) [D]^2.
\end{equation}
This, in turn, says that the parabolic Chern classes of $E$ are zero.

\end{proposition}

If $D$ is smooth and $X$ is a surface, this follows from \cite{B2}. 
The same proof works here
with the only variation that the parabolic structure is defined on the
normalization $\wt D$ of $D$. The details of the proof are given in
section  \ref{proofbiquard}.

\begin{lemma}
\label{mumfordunstable}
With the same notation as above, assume that $D^2>0$ and that not all
weights $\alpha_i$ are equal. Then $E$ is Mumford unstable.
\end{lemma}

\begin{proof}
Using (\ref{chernchar}) we calculate the discriminant of $E$ 
\begin{eqnarray*}
\lefteqn{\Delta= \Big( -\ch_2(E) + \frac{1}{2r}(\ch_1(E))^2 \Big)
H^{n-2} =}\\
& & -\frac{1}{2}\Big(\sum \alpha_i^2 - \frac{1}{r}(\sum \alpha_i)^2\Big)
D^2 H^{n-2}<0,
\end{eqnarray*}
then, by Bogomolov inequality (\cite[thm 7.3.1]{H-L}), 
$E$ is Mumford unstable.

\end{proof}

{}From now on we will assume that $\NS(X)=\ZZ L$.

\begin{proposition}
\label{parabolicunstable}
Assume (\ref{assumption}).
Let $E$ be a 
parabolic vector bundle with vanishing parabolic Chern
classes and rank $r$. 
If $E$ is parabolic stable, then the parabolic structure is trivial
(all the weights are equal).
\end{proposition}

\begin{proof}
Assume that not all weights are equal.
By lemma \ref{mumfordunstable}, $E$ is Mumford unstable. We
will prove that $E$ is not parabolic stable by showing that at least one
of the subsheaves of the Harder-Narasimhan filtration contradicts the
parabolic stability of $E$. Since the parabolic degree of $E$ is zero, we
have to prove that the parabolic degree of one of the subbundles is
non-negative.

First we assume that the Harder-Narasimhan filtration of $E$ has only
one term, i.e. there is a short exact sequence
\begin{equation}
0 \to E' \to E \to E'' \to 0
\label{shortexact}
\end{equation}
such that both $E'$ and $E''$ are Mumford semistable torsion-free
sheaves. 
The objective is
to show that the parabolic degree of $E'$ is non-negative.

Let $s=(1/r)\sum \alpha_i$, and $\alpha^0_i = \alpha_i - s$. By the
formula for the first Chern class of $E$ we have $c_1(E)\equiv -rsdL$,
where ``$\equiv$'' means numerical equivalence, and $d L\equiv [D]$. 
Let $a'L\equiv 
c_1(E')$, $r'=\rk(E')$, and analogously for $a''$
and $r''$.

Both $E'$ and $E''$ are Mumford semistable, then by Bogomolov
inequality
$$
0 \leq \Big( c_2(E')-\frac{r'-1}{2r'}{a'}^2L^2 +
c_2(E'')-\frac{r''-1}{2r''}{a''}^2L^2 \Big) H^{n-2}.
$$
Using $c_2(E)=c_2(E')+c_2(E'')+a'a''L^2$ and the formula for $c_2(E)$
in terms of $\Gamma$ this inequality becomes
$$
0 \leq \Big((r^2-r)s^2 d^2-\sum(\alpha_i^0)^2 d^2- 2 a'a'' - 
\frac{r'-1}{r'}{a'}^2 - \frac{r''-1}{r''}{a''}^2\Big)\frac{L^2}{2}
H^{n-2}
$$
Using the formula for $c_1(E)$ in terms of $\Gamma$ we have 
$a''=-rsd-a'$. Substituting this into the inequality and simplifying
we obtain
$$
0 \leq \Big(\frac{r r'}{r''}(\frac{a'}{r'}+sd)^2 - \sum(\alpha^0_i)^2
d^2 \Big)
\frac{L^2}{2} H^{n-2}
$$
Note that because $E'\subset E$ is a Harder-Narasimhan filtration we
have
$$
\frac{a'}{r'}+sd=\frac{a'}{r'}-\frac{a}{r}>0,
$$
and then we have
\begin{equation}
\label{ine1}
\frac{a'}{r'}+sd-\sqrt{\frac{r''}{r'r}\sum(\alpha^0_i)^2d^2}\geq 0.
\end{equation}
There is an induced parabolic structure on $E'$. There is a subset 
$I'\subset I=\{1,\ldots,r\}$ with $r'$ elements such that the 
weights of the parabolic
structure on $E'$ induced by $E$ are $\alpha_i$ for $i \in I'$.

\textbf{Claim.} For any subset $I'\subset I$ of cardinality $r'$,
\begin{equation}
\label{ine2}
\frac{\sum_{i\in I'} \alpha^0_i}{r'} \geq 
-\sqrt{\frac{r''}{r'r}\sum_{i\in I}(\alpha^0_i)^2}.
\end{equation}
Proving this is an easy calculus exercise. Use the method of Lagrange
multipliers to minimize $\sum_{i\in I'} \alpha^0_i$ subject
to the conditions $\sum_{i\in I} \alpha^0_i=0$ and 
$\sum_{i\in I}(\alpha^0_i)^2= R$ where $R$ is some constant (note that
we are minimizing a linear function restricted to a sphere).

Then combining inequalities (\ref{ine1}) and (\ref{ine2}) we get 
$$
0 \leq \frac{a'+\sum_{i\in I'} \alpha^0_i}{r'}+sd =
\frac{a'+\sum_{i\in I'} \alpha_i}{r'},
$$
but this is the parabolic degree of $E'$, so $E$ cannot be 
parabolic stable.

Now we assume that the Harder-Narasimhan filtration has length 2,
i.e. we have
$$
E_1 \subset E_2 \subset E
$$
Then $\rk(E)=3$ (recall that we assume $r\leq 3$), $\rk(E_i)=i$, and
$E_2/E_1$ and $E/E_2$ are torsion free sheaves of rank one. Let
$a_1=c_1(E_1)$, $a_2=c_1(E_2/E_1)$, $a_3=c_1(E/E_1)$ (since
$NS(X)\isom \ZZ L$, we can think of the first Chern class as an
integer number). By the
definition of the Harder-Narasimhan filtration we have 
\begin{equation}
\label{HNineq}
a_1>a_2>a_3
\end{equation}
Using the formula (\ref{chern}) for $c_1(E)$ we have $-\sum \alpha_i d =
c_1(E)=\sum a_i$, where $d$ is the degree of $D$. Then there exist
(rational) numbers $x$, $y$ such that
\begin{equation}
\label{alphas1}
\alpha_1 d = -\frac{c_1(E)}{3}+x,\quad 
\alpha_2 d = -\frac{c_1(E)}{3}+y,\quad
\alpha_3 d = -\frac{c_1(E)}{3}-x-y.
\end{equation}
Using the formula for $\ch_2(E)$ we have
\begin{equation}
\label{alphas2}
\sum \alpha_i^2 d^2 L^2= 2 \ch_2(E) = \Big(-l + \sum a_i^2\Big)L^2,
\end{equation}
where $l L^2=2 c_2(E_2/E_1) +2 c_2(E/E_1)$, hence $l \geq 0$. 
Combining this with
(\ref{alphas1}) we obtain
$$
y=-\frac{x}{2} \pm \sqrt{\frac{3}{4}(x_m^2-x_{}^2)},\;
\text{where}\;\;
x_m^2=\frac{2}{3}\Big(-l+\sum a_i^2\Big)-
\frac{2}{9}\Big(\sum a_i\Big)^2.
$$
The number $y$ must be real, then $-x_m\leq x \leq x_m$.

Now we calculate the parabolic degrees of $E_1$ and $E_2$ as
functions of $x$ (we take the polarization $(1/L^2) L$). Using
$\pardeg E_1=a_1+\alpha_1 d$ and $\pardeg
E_2=a_1+\alpha_1 d+ a_2+\alpha_2 d$ we obtain
\begin{eqnarray}
\label{pardeg1}
\pardeg E_1 & = & x+\frac{2a_1-a_2-a_3}{3} \\
\label{pardeg2}
\pardeg E_2 & = & \frac{a_1+a_2-2a_3}{3}
+\frac{x}{2}\pm \sqrt{\frac{3}{4}(x^2_m-x^2_{})} 
\end{eqnarray}

Note that if we fix $x$, $\pardeg E_1$ is fixed, but $\pardeg E_2$ could
take two values, hence the two signs in the formula. We want to show
that for any value of $x$ (with $x^2_{}\leq x^2_m$), at least one of
these is non-negative.

\begin{figure}[ht]
\centerline{\epsfig{file=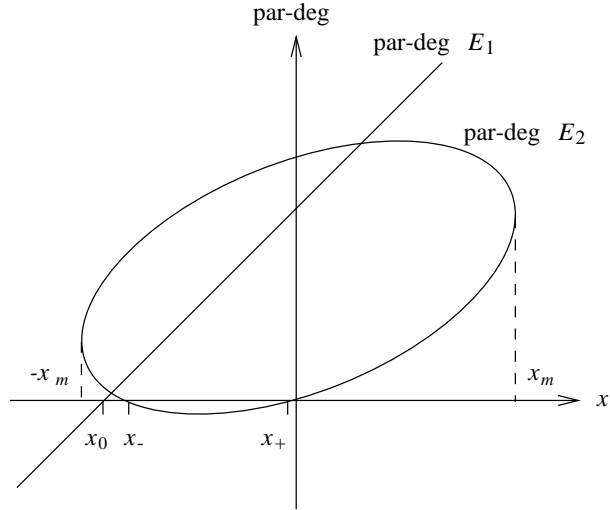}}
\caption{$x_0=\frac{-2a_1+a_2+a_3}{3}$,
$x_\mp=\frac{-a_1-a_2+2a_3}{6}\mp \frac{1}{2}\sqrt{(a_1-a_2)^2-2l}$
}
\end{figure}

Equation (\ref{pardeg2}) defines a conic with
coordinates $(x,\pardeg E_2)$ in the plane $\RR^2$.  
This conic intersects the axis
$\pardeg E_2=0$ in the two points $x_-$ and $x_+$
$$
x_\pm=\frac{-a_1-a_2+2a_3}{6}\pm \frac{1}{2}\sqrt{(a_1-a_2)^2-2l}
$$
For $x=x_m$, $\pardeg E_2= \frac{1}{3}
(a_1+a_2-2a_3)+\frac{1}{2}x_m>0$ (by equation (\ref{HNineq})). 
Assume $\pardeg E_2<0$. Then 
$$
x_-<x<x_+,
$$
but if $x$ is in this interval, using $l\geq 0$ we have 
$$
x>x_->\frac{-a_1-a_2+2a_3}{6}- \frac{1}{2}\sqrt{(a_1-a_2)^2}
=\frac{-2a_1+a_2+a_3}{3}
$$
(for the equality we used equation (\ref{HNineq}), to get the correct
sign of the square root), and then equation (\ref{pardeg1}) 
implies that $\pardeg E_1>0$.
Then either $\pardeg E_1$ or $\pardeg E_2$ is non-negative.

\end{proof}

\begin{corollary}
\label{corollary}
Assume (\ref{assumption}).
If $\rho:\pi_1(X\minus D)\to U(r)$ is an irreducible
representation then $\rho(\lambda)$ is a 
multiple of the identity.
\end{corollary}

\begin{proof}
Let $E$ be the corresponding stable parabolic bundle given by
proposition \ref{biquard}. By proposition
\ref{parabolicunstable} the parabolic structure is trivial (i.e. all
weights are equal), and then $\rho(\lambda)$ is a multiple of
identity.

\end{proof}

\begin{theorem}
\label{irred1}
Assume (\ref{assumption}).
If $\pi_1(X)$ has no irreducible $PU(r)$ representations 
then $\pi_1(X \minus D)$ has no 
irreducible $U(r)$ representations.
\end{theorem}

\begin{proof} Let $\rho:\pi_1(X \minus D)\to U(r)$ be a representation,
and $E$ the associated parabolic bundle. If $\rho$ is irreducible
then, by corollary \ref{corollary}, $\rho(\lambda)$ is a multiple of
identity.  
Then the 
representation $\rho:\pi_1(X\minus D)\to U(r)$ induces a
representation $\overline\rho:\pi_1(X)\to PU(r)$. If $\rho$ is
irreducible then $\overline\rho$ is also irreducible, and we get a
contradiction.

\end{proof}

\section{Non-unitary representations of fundamental groups}

One can also ask about non-unitary representations. In this section we
will assume that $X$ is a surface and $D$ is smooth. 
We only deal with the
case of $GL(2)$ representations.

Recall (\cite[d\'ef 1.1]{B3}) that a parabolic Higgs bundle is a 
parabolic bundle $E$
together with a section $\phi\in H^0(\Omega^1\langle \log D \rangle 
\otimes ParEnd E)$ with $\phi_\wedge \phi=0$, 
where $ParEnd E$ is the sheaf of
parabolic endomorphisms of $E$. The residue $\Res_D \phi$ respects
the filtration $F_\bullet$ of $E|_D$, and we require that its conjugation
class is constant on each quotient $F_j/F_{j+1}$.

The following result follows from Biquard's theorem \cite[thm 11.4]{B3}.

\begin{proposition}
\label{frombiquard}
Let $D$ be a smooth divisor with $D^2\neq 0$. 
Given an irreducible local system of rank $r$ on $X-D$, 
there is a rank $r$
stable parabolic Higgs bundle $(E,\phi)$ with vanishing parabolic
Chern classes. 

The eigenvalues of the residue of $\phi$ are zero.
The slopes of the bundles $F_j/F_{j+1}$ are given by 
\begin{equation}
\label{degreegr}
\frac{\deg(F_j/F_{j+1})}{\rk(F_j/F_{j+1})}=-\mu_j D^2,
\end{equation}
where $\mu_j$ are the parabolic weights of the parabolic Higgs bundle.

\end{proposition}

\begin{proof}
Given an irreducible local system on $X-D$, consider its Deligne
extension $(E,\nabla)$ (\cite{Ka}). By construction, the real part of
the eigenvalues
of the residue of $\nabla$ are non-negative and less than 1. Define
a parabolic structure on $E$, setting the parabolic weights equal to
the real parts of the eigenvalues of the residue of $\nabla$. This
integrable logarithmic connection satisfies the hypothesis of
\cite[thm 11.4]{B3}, and then we obtain a parabolic Higgs bundle.
It is stable because of the irreducibility hypothesis and because it 
has a Hermite-Einstein metric.

The vanishing of the 
eigenvalues of the residue of $\phi$ is given in \cite[lemme 7.1]{B3}.
By \cite[lemme 11.2]{B3}, the eigenvalues of the residue of $\nabla$
are real, and by \cite[prop 10.1]{B3} they are equal to the parabolic
weights of the parabolic Higgs bundle, and applying again \cite[lemme
11.2]{B3} we obtain the formula for the slopes.

\end{proof}

We use this result to prove the following theorem.

\begin{theorem}
\label{irred2}
Assume that $X$ is a projective surface, $D$ is a smooth divisor and  
$\NS(X)\isom \ZZ L$. If $\pi_1(X)$ has no
irreducible $PGL(2)$ representations, then $\pi_1(X \minus D)$ has no
irreducible $GL(2)$ representations.
\end{theorem}

\begin{proof} 
Assume that there is an irreducible representation  $\rho:\pi_1
(X \minus D) \to GL(2)$.
Let $(E,\phi)$ be the stable parabolic Higgs bundle associated by
proposition \ref{frombiquard}.

We will use the following fact: 
\begin{equation}
\label{logpole} 
H^0(\Omega^1 \langle \log D\rangle)=H^0(\Omega^1).
\end{equation} 
The proof was shown to us by I. Biswas. Consider the sequence 
\begin{equation}
\label{shortresidue}
0 \to \Omega^1 \to \Omega^1 \langle \log D\rangle \to \SO_D \to 0,
\end{equation} 
where the second map is the residue, and note that the image of the 
constant function
1 under the co-boundary map is represented by the \v{C}ech cocycle 
$dz_i/z_i -dz_j /z_j=d \log (z_i/z_j)$, where $z_i$ is a local equation
of $D$ on an open set $U_i$, and this is
(up to a non-zero constant) the Chern class
of the line bundle $\SO(D)$ (cf. \cite[prop 12]{A}).

If $M$ is a nontrivial line bundle on $X$ with $\deg M\leq 0$, we
still have 
\begin{equation}
\label{noresidue}
H^0(M\otimes \Omega^1 \langle \log D\rangle)=H^0(M\otimes \Omega^1).  
\end{equation} 
To see this, note that $H^1(M(-D))=0$ by Kodaira vanishing theorem,
hence $H^0(M|_D)=H^0(M)=0$, and then use the sequence (\ref{shortresidue})
tensored by $M$.

Consider first the case when the parabolic structure of $E$ is non-trivial
(weights are different).  We
already know that the eigenvalues of $\Res_D \phi$ are zero. Since the
residue of the Higgs field also preserves the parabolic filtration, it
yields a map of line bundles $F_1/F_2 \to F_2$. But by
(\ref{degreegr}) we have
$$ 
\deg(F_1/F_2)=-\mu_1 D^2 > -\mu_2D^2 = \deg(F_2)  
$$ 
so this map has to be zero. In other words, $\Res_D \phi=0$, and then 
$(E,\phi)$ defines a  Higgs pair on $X$. Using Bogomolov inequality for
Higgs bundles (cf. \cite[prop 3.4 and thm 1]{Si}) 
we see that $(E,\phi)$ is not stable as a 
Higgs pair. Then there is a short exact sequence
$$
0 \to M \to E \to M' \otimes I_Z \to 0
$$
where $I_Z$ is the ideal sheaf of a zero-dimensional subscheme and 
$M$ and $M'$ are line bundles with $\deg M \geq \deg M'$. A
calculation similar to the proof of proposition \ref{parabolicunstable}
(for filtrations of length 1), using formulae (\ref{chern}) and the
fact that $c_2(I_Z)=l(Z)\geq 0$, shows that $E$ is not parabolic Higgs
stable (since $M$ is $\phi$-invariant), contradicting the 
hypothesis and then finishing the proof in this case.

The case when the parabolic structure is trivial is more subtle (it is to
cover this case that we need the assumption on the tangent bundle).
Since all the weights are equal, $E$ is
Mumford stable iff it is parabolic stable. But if $E$ is parabolic
stable then there is a corresponding $U(2)$ irreducible representation
$\rho^{}_{U(2)}$ of $\pi_1(X\minus D)$, and since the weights are
equal, $\rho^{}_{U(2)}(\lambda)$ is a multiple of the identity. 
Then
there is an induced irreducible $PU(2)$ representation of $\pi_1(X)$, in
contradiction with the hypothesis.

Then $E$ is not Mumford stable, and thus there is a sequence as before
$$
0 \to M \to E \to M' \otimes I_Z \to 0 
$$
where $M$ and $M'$ are line bundles, $I_Z$ is the ideal
sheaf of a zero-dimensional subscheme $Z$, and 
$\deg M \ge \deg M'$, but now $M$ might not be $\phi$-invariant, so
$M$ doesn't contradict parabolic Higgs stability. 
The Higgs field $\phi$
defines a map $\phi':M \to M' \otimes \Omega^1 \langle
\log D\rangle$.  Regard
this as a section of $H^0({M}^{-1} \otimes M' \otimes
\Omega^1 \langle \log D\rangle)$. By (\ref{noresidue}),
this is a section of $H^0({M}^{-1} \otimes M' \otimes
\Omega^1)$. In other words, $\Res \phi'=0$. On the other hand, $\phi'$
vanishes on $Z$, since it factors through $M'\otimes I_Z\otimes
\Omega^1 \langle \log D \rangle$. Putting both facts together, $\phi'$
can be seen as a section of $H^0(M^{-1}\otimes M' \otimes I_Z \otimes
\Omega^1)$. If the tangent bundle is globally generated we see that this
group is zero (and thus $M$, being invariant under $\phi$,
contradicts stability of the parabolic
Higgs pair) unless $M= M'$ and $Z$ is empty. 

Suppose therefore that we have a non-trivial extension 
\begin{equation}
\label{ext1}
0 \to M \to E \to M \to 0 
\end{equation}

We now produce a ``model'' of the extension (\ref{ext1}), which is adapted
to the divisor $D$.  
Since the tangent bundle is generated by global sections, 
the anticanonical line bundle $K^{-1}$ has non-negative degree. Then,
since the Neron-Severi group has rank one, $K^{-1} \otimes \SO(D)$ is
ample, and by Kodaira vanishing theorem $H^1(\SO(D))=0$.

Let $\SN_D=\SO(D)|_D$ denote the normal bundle of
$D$. From the sequence 
\begin{equation} 
\label{extensiond}
0 \to \SO \to \SO(D) \to \SN_D \to 0 
\end{equation} 
it follows that the
map $H^0(\SN_D) \to H^1(\SO)$ is onto. 
Note that this map is obtained from the natural multiplication 
map $\Hom(\SO,\SN_D)\times 
\Ext^1(\SN_D,\SO) \to \Ext^1(\SO,\SO)$, using the extension class of
(\ref{extensiond}).
Then, given a section $\sigma$ of $\SN_D$,
the corresponding extension $E_{\sigma}$ is given as the kernel in the
sequence
\begin{equation}
\label{ext2}
0 \to E_{\sigma} \to \SO(D) \oplus \SO \to \SN_D \to 0 
\end{equation}
the map on the right being $(u,v)\mapsto u|^{}_D- \sigma v|^{}_D$.
From (\ref{ext1}) we have
$$
0 \to \SO \to E\otimes M^{-1} \to \SO \to 0,
$$
and since $H^0(\SN_D) \to H^1(\SO)$ is surjective, there is a section
$\sigma$ such that $E=E_\sigma \otimes M$. From now on we will use this
``model'' of $E$.

Note the obvious inclusion $M \oplus M (-D) \hookrightarrow
E$. Consider the composition 
$$
M \oplus M (-D) \inj E \stackrel{\phi}\to
E \otimes \Omega^1\langle \log D\rangle \inj \{M (D) \oplus M \}
\otimes \Omega^1 \langle \log D\rangle,
$$
where the last map comes from the left map on (\ref{ext2}). We denote
this composition by $\hat\phi$. One can now represent this map as a matrix 
\begin{equation}
\label{higgs}
\left(\matrix
a &b \\
c &d\\
\endmatrix\right)
\end{equation}
where $a,d$ are sections of $\SO(D) \otimes \Omega^1 \langle \log
D\rangle$, 
$c$ is a section
of $\Omega^1 \langle \log D\rangle$, and $b$ of $\SO(2D) \otimes \Omega^1 
\langle \log D\rangle$. We will
use this representation to show that the residue of the Higgs field is
zero (note that $\Res \phi = \Res \hat\phi$).

By (\ref{logpole}), $c \in
H^0(\Omega^1)$. Consider a local holomorphic section $(f,0)$ of 
$M \oplus M (-D)$. The residue of $\hat\phi$ sends 
$(f|_D,0)$ to 
$$
(\Res(a)f|_D,\Res(c)f|_D)=(\Res(a)f|_D,0),
$$ 
hence
$\Res(a)$ is an eigenvalue of $\Res(\hat\phi)$, and hence of
$\Res(\phi)$, but those eigenvalues are zero, then $\Res(a)=0$ and
$a\in H^0(\SO(D)\otimes \Omega^1)$.

Since the eigenvalues of $\hat\phi$ are zero, $0=\tr \Res(\hat\phi)=
\Res(a)+\Res(d)=\Res(d)$, and then $d\in H^0(\SO(D)\otimes \Omega^1)$.

Finally we will impose the condition that the map $\hat\phi$ comes
from the Higgs field on $E$, i.e. $\hat\phi$ extends to $E$. The exact
sequence (\ref{ext2}) shows (after tensoring with $M$) that a 
local section of $E$ is the same thing as a local section $(f',g')$  
of $M(D) \oplus M$, with the condition that $f'|_D-
\sigma g'|_D=0$. We will write local sections of $M(D) \oplus
M$ as $((1/z)f,(1/z)g)$, where $f$ and $g$ are respectively
holomorphic local sections of $M$ and $M (-D)$ and $z$
is a local equation of $D$. Take a local holomorphic section $g$ of
$M (-D)$ that doesn't vanish identically on $D$, and let $f$ be
a local holomorphic section of $M$ such that $f|_D-
\sigma g|_D=0$, so that $((1/z)f,(1/z)g)$ is in fact a local section
of $E$. The matrix (\ref{higgs}) acting on this gives 
$$
( \frac{1}{z} af+  \frac{1}{z} bg, \frac{1}{z} cf + \frac{1}{z} dg).
$$
This should be a local section of $\{M (D) \oplus M \}
\otimes \Omega^1 \langle \log D\rangle$. 
We have seen that $a$ is in fact a section of $\SO(D)\otimes
\Omega^1$, then $(1/z)a f$ is a local section of $M
(D)\otimes \Omega^1\langle \log D\rangle$. We need then $(1/z)b g$ to be 
also a local
section of $M(D)\otimes \Omega^1\langle \log D\rangle$. 
Recall that $b$ is a
section of $\SO(2D)\otimes \Omega^1\langle 
\log D\rangle$ and $g$ a local section of
$M(-D)$ (that doesn't vanish identically on $D$). Then for 
$(1/z) b g$ to be a local section of $M(D)\otimes
\Omega^1\langle \log D\rangle$ we need the residue of $b$ to be zero, 
i.e. $b\in
H^0(\SO(2D)\otimes \Omega^1)$. Then we obtain that the residue of
(\ref{higgs}) (and hence of $\phi$) is zero.

Finally, if the residue of $\phi$ is zero, since we have assumed that
the weights of the parabolic structure of $E$ are equal, we obtain
that $\rho(\lambda)$ is a multiple of the identity, then there is an
induced $PGL(2)$ representation, and since this is reducible by 
hypothesis, we conclude that $\rho$ is also reducible.

\end{proof}

\section{Proof of proposition \ref{biquard}}
\label{proofbiquard}

We will follow closely \cite[(V.8.5)]{Ko} and \cite[(4.1)]{B2}. This
proof works for $X$ a K\"ahler manifold, not only for projective
varieties. We will denote by $\omega$ the K\"ahler form.

The representation $\rho$ gives a local system on $X\minus D$. Let
$(E,\nabla)$ be the Deligne extension.
Let $p$ be a point of $D$. There are local coordinates $z_1,\ldots,z_n$
on $X$ and a local holomorphic trivialization of $E$, 
where $D$ is defined by
$t=z_1\cdots z_l$, and the Deligne extension is (\cite{Ka})
$$
\nabla= d+ \sum_{j=1}^{l} \frac{dz_j}{z_j} \Gamma_j
$$
where $\Gamma_j=\diag(\alpha^{(j)}_1,\ldots,\alpha^{(j)}_r)$, and 
$\alpha^{(j)}_1,\ldots,\alpha^{(j)}_r$ is a permutation of
$\alpha_1,\ldots,\alpha_r$.

Fix once and for all a point $x_0\in X$ and a Hermitian metric on the
fiber $E|_{x_0}$. Since the connection is flat and the holonomy is
unitary, using this fixed metric and parallel transport we can define a
Hermitian metric on $E$ (degenerate on $D$). 
In the previous local coordinates
$$
h=\left(
\begin{array}{ccc}
|z_1|^{\alpha^{(1)}_1} \dots |z_l|^{\alpha^{(l)}_1} & & 0 \\
 & \ddots & \\
0 & & |z_1|^{\alpha^{(1)}_r} \dots |z_l|^{\alpha^{(l)}_r}
\end{array}
\right)
h_0
$$
where $h_0$ is a fixed constant matrix (depending on the metric
chosen on $E|_{x_0}$). This is an ``adapted metric'', in the sense of
\cite[d\'ef 2.3]{B1}. By direct computation we can check that $\nabla$
is the Chern connection associated to the metric $h$.

Now we define the parabolic structure (\cite[thm 2.1]{B2}). Since the
local equation for $D$ is $t=z_1\cdots z_l=0$, where $l$ can be
different from 1 in general,
we need to go to the normalization $\wt D$ of $D$ to separate the different
branches. Let
$\tilde p \in \wt D$ be a point mapping to $p\in D$, such that a small
neighborhood $\wt V$ of $\wt D$ maps to the subset of $D$ defined
by the equation $z_k=0$, with $1\leq k\leq l$. We define a filtration
of $E|_{\wt V}=\pi^*_{\wt V} E$
$$
E|_{\wt V}=F_1|_{\wt V} \supset F_2|_{\wt V} \supset \ldots \supset
F_r|_{\wt V} \supset 0
$$
by the property that if $s(z)$ is a local section of $E$, then 
$$
s|_{\wt V} \in F_i|_{\wt V}\minus F_{i+1}|_{\wt V} 
\quad\Leftrightarrow\quad 
\|s(z)\|_h \sim |z_k|^{\alpha_i}
$$
i.e. if $\|s(z)\|_h=|z_k|^{\alpha_i} f$, with $f$ a continuous
function that doesn't vanish identically on the branch $z_k=0$.
This local construction defines a filtration on $E|_{\wt D}$ and hence
a parabolic structure on $E$.

Now we will prove that if the
representation $\rho$ is irreducible, this parabolic structure is
stable (it will then follow that if $\rho$ is not irreducible, then
the parabolic structure will still be polystable: decompose $\rho$ in
irreducible representations, and then the parabolic bundle will be a
direct sum of stable parabolic bundles, i.e. a polystable parabolic bundle).

Let $E'$ be a saturated coherent torsion-free subsheaf of $E$ of rank
$r'$. There is a naturally induced parabolic structure 
on $E'$ whose weights are a subset of the weights of $E$ (to check
stability it is enough to look at parabolic subsheaves with this parabolic
structure, since they have maximal parabolic degree).
It is defined by considering the Hermitian metric $h'$ on $E'$
obtained by restriction of
the metric $h$ on $E$. The metric $h'$ is degenerate on the points
where $E'$ is not locally free or where $i:E'\to E$ is not an
injection of vector bundles. Let $S$ be the union of all such
points. We define a parabolic structure as before, but only on
$X\minus S$ (i.e. the filtration of $\pi^*E'=E'|_{\wt D}$ is only
defined on $\wt D \minus \pi^{-1}(S)$. This is enough for our purpose:
to calculate the parabolic degree).
We denote by $\sum_{i\in I'}\alpha_i$ the sum of the parabolic weights
of $E'$ (with repetitions), where $I'\subset\{1,\ldots,r\}$.

\begin{lemma}
\begin{equation}
\label{curvature}
\frac{\sqrt{-1}}{2\pi}\int_{X\minus (D\cup W)} \tr \delb({h'}^{-1}\del h')
{}_\wedge \omega^{n-1} = \pardeg E'
\end{equation}
\end{lemma}

\begin{proof}
Since $h'$ (and also the associated Chern connection) is singular on
$D$, we cannot directly apply the Chern-Weil formula. We will modify
the metric to make it smooth, and this will produce the second summand.

First note that 
\begin{eqnarray}
\label{determinant}
\tr \delb({h'}^{-1}\del h')=\delb\del \log (\det h'). 
\end{eqnarray}
This can be proved by first showing that both sides are
invariant under change of holomorphic trivialization, and then
computing in a holomorphic trivialization where $h=\id + O(|z|^2)$ (\cite[III
lemma 2.3]{W}). Note that $\det h'$ is a Hermitian metric on $\det
E'$, singular on $D$ (in fact $\det h' \sim |z_1\cdots z_l|^{2
(\sum_{i\in I'}\alpha_i)}$).

Following \cite[(4.1)]{B2}, choose a section $t\in H^0(\SO_X(D))$ with
a zero of order one along $D$.
Take any Hermitian metric on the line bundle
$\SO_X(D)$, and then $\|t\|$ is a smooth function on $X$, vanishing on
$D$. Near $D$, $t=z_1\cdots z_l$ (using an appropriate trivialization
of $\SO_X(D)$), and then 
$$
\frac{\det h'}{\|t\|^{2(\sum_{i\in I'}\alpha_i)}}
$$
defines a smooth metric on $(\det E')|_{X\minus S}$, 
and we can apply the Chern-Weil
formula as in \cite[(V.8.5) formula (**)]{Ko}.
The integral (\ref{curvature}) can be written, using 
(\ref{determinant})
\begin{eqnarray*}
& & \frac{\sqrt{-1}}{2\pi}\int_{X\minus (D\cup S)} 
\Big(
\delb\del \log \frac{\det h'}{\|t\|^{2(\sum_{i\in I'}\alpha_i)}}
+(\sum_{i\in I'}\alpha_i)
 \delb\del \log \|t\|^2
\Big)
{}_\wedge \omega^{n-1}    \\
& &  = \deg E' +(\sum_{i\in I'}\alpha_i)\deg D 
\end{eqnarray*}
where the second integral is given by the Poincar\'e-Lelong formula.
This finishes the proof of the lemma.

\end{proof}

If we take $E'=E$, then we obtain $\pardeg E=0$, since the
connection is flat on $X \minus D$. Using this lemma, to prove that
$E$ is stable we have to show that for any proper saturated
torsion-free subsheaf $E'$ of
$E$,
\begin{eqnarray}
\label{stability}
\frac{\sqrt{-1}}{2\pi}\int_{X\minus (D\cup S)} \tr \delb({h'}^{-1}\del h')
{}_\wedge \omega^{n-1} <0.
\end{eqnarray}
First note that the Chern connection associated to $h'$ is equal to
$\nabla'=\pi_{E'}\circ \nabla$, where $\pi_{E'}$ is the orthogonal
projection (with respect to $h$) on the subbundle $E'$. Then 
$$
\delb({h'}^{-1}\del h)=\Theta_{\nabla'}=\Theta_{\nabla}|_{E'}-A_\wedge
A^*,
$$
where $\Theta_{\nabla'}$ is the curvature of $\nabla'$, $A$ is the
second fundamental form and the second
equality if given by the Gauss-Codazzi formula (see \cite{Ko}).
But $\Theta_\nabla=0$ on $X\minus D$, since $\nabla$ is flat, and
then the integral (\ref{stability}) can be written as
$$
\frac{\sqrt{-1}}{2\pi}\int_{X\minus (D\cup S)} - \tr A_\wedge A^* 
{}_\wedge\omega^{n-1}
$$

For any $A$, this integral is non-positive. It is zero only if 
the second fundamental form $A$ is
identically zero, but this would imply that there is a holomorphic
splitting $E|_{X\minus (D\cup S)}=E'|_{X\minus (D\cup S)} \oplus
{E'}^{\perp}|_{X\minus (D\cup S)}$. The same argument in 
\cite[(V.8.5) p. 183]{Ko}
shows that this extends to a holomorphic splitting on $X-D$: since we
have a splitting on $X\minus (D\cup S)$, the holonomy group 
$G$ of $(E,h)$ on
$X \minus (D\cup S)$ is in $U(r')\times U(r-r')$, and since $S$ is a
closed subvariety (of $\codim S\geq 2$), the holonomy group on $X\minus
D$ is contained in the closure of $G$, and this is still in
$U(r')\times U(r-r')$. 

This contradicts the irreducibility
of the representation $\rho$. Then this integral has to be negative, and
(\ref{stability}) is proved. This finishes the proof of proposition 
\ref{biquard}.

\bigskip
\bigskip

\noindent\textbf{Acknowledgments.}
We would like to thank I. Biswas for discussions. We would also like 
to thank O. Biquard for explaining us some points of his work.
T.G. was supported by a postdoctoral fellowship of 
Ministerio de Educaci\'on y Cultura (Spain).


\begin{thebibliography}{EMG} 

\bibitem[A]{A}{Atiyah, M.F.:}
Complex analytic connections in fibre bundles,
Trans. Amer. Math. Soc. \textbf{85}, 181-207 (1957)

\bibitem[B1]{B1}{Biquard, O.:}
Fibr\'es paraboliques stables et connexions singuli\`eres
plates, 
Bull. Soc. Math. France \textbf{119}, 231--257 (1991)

\bibitem[B2]{B2}{Biquard, O.:}
Sur les fibr\'es paraboliques sur une surface complexe, 
J. London Math. Soc. \textbf{53}, 302--316 (1996)

\bibitem[B3]{B3}{Biquard, O.:}
Fibr\'es de Higgs et connexions int\'egrables: le cas
logarithmique (diviseur lisse), 
Ann. scient. \'Ec. Norm. Sup., 4 s\'erie, \textbf{30},
41--96 (1997)

\bibitem[D]{D}{Deligne, P.:}
Equations differentielles \`a points singuliers r\'eguliers, 
Lecture Notes in Mathematics \textbf{163}, Springer Verlag, 1970

\bibitem[E-V]{E-V}{Esnault, H., Viehweg, E.:}
Logarithmic de Rham complexes and vanishing theorems, 
Invent. math. \textbf{86}, 161--194 (1986)

\bibitem[H-L]{H-L}{Huybrechts, D., Lehn, M.:}
The geometry of moduli spaces of sheaves, 
Aspects of Mathematics E31, Vieweg, Braunschweig/Wiesbaden 1997.

\bibitem[Ka]{Ka}{Katz, N.:}
An overview of Deligne's work on Hilbert's twenty-first
problem, 
in \textit{Mathematical developments arising from Hilbert problems, }  
Proc. Sympos. Pure Math., Vol. XXVIII, Amer. Math. Soc.,  537--557 (1976)

\bibitem[Ko]{Ko}{Kobayashi, S.:}
Differential geometry of holomorphic vector bundles, Princeton,
Princeton University Press, 1987

\bibitem[K-M]{K-M}{Kronheimer, P., Mrowka, M.:}
Gauge theory for embedded surfaces, I, 
Topology \textbf{32}, 773--826 (1993)



\bibitem[M-Y]{M-Y}{Maruyama, M., Yokogawa, K.:}
Moduli of parabolic stable sheaves, 
Math. Ann. \textbf{293}, 77-99 (1992)

\bibitem[No]{No}{Nori, M.:}
Zariski's conjecture and related
problems,  Ann. scient. \'Ec. Norm. Sup., 4 s\'erie, \textbf{16}, 305-344
(1983)

\bibitem[Si]{Si}{Simpson, C.:}
Constructing variations of Hodge
structure using Yang-Mills theory and applications to uniformization, 
Journal of the AMS \textbf{1}, 867-918 (1988)

\bibitem[W]{W}{Wells, R.O.:}
Differential analysis on complex manifolds, 
Grad. Texts in Math. \textbf{65}, Springer Verlag, 1980



\end{thebibliography}
\end{document}